\documentclass[a4paper]{scrartcl}

\usepackage{amsfonts, amssymb, amsthm}

\usepackage{graphicx, xcolor}
\usepackage{color}

\usepackage[english]{babel}
\usepackage{mathtools}
\usepackage{booktabs}
\usepackage{hyperref}
\mathtoolsset{showonlyrefs,showmanualtags}

\newcommand{\onelinewidth}{0.84\linewidth}

\newcommand{\cL}{\mathcal{L}}

\newcommand{\cH}{\mathcal{H}}

\newcommand{\der}[2]{\frac{\partial{#1}}{\partial{#2}}}

\newcommand{\dd}{\mathrm{d}}

\newcommand{\R}{\mathbb{R}}

\newtheorem{definition}{Definition}
\newtheorem{proposition}{Proposition}

\theoremstyle{definition}

\title{New directions for contact integrators}

\usepackage{authblk} 
\usepackage{orcidlink}

\author[1]{Alessandro Bravetti \orcidlink{0000-0001-5348-9215}}
\author[2]{Marcello Seri \orcidlink{0000-0002-8563-5894}}
\author[2]{Federico Zadra \orcidlink{0000-0001-7565-3792}}

\affil[1]{Instituto de Investigaciones en Matemáticas Aplicadas y en Sistemas (IIMAS--UNAM), Mexico City, Mexico}
\affil[2]{Bernoulli Institute for Mathematics, Computer Science and Artificial Intelligence, Groningen, The Netherlands}

\date{}

\begin{document}

\maketitle

\abstract{
    Contact integrators are a family of geometric numerical schemes which guarantee the conservation of the contact structure.
    In this work we review the construction of both the variational and Hamiltonian versions of these methods.
    We illustrate some of the advantages of geometric integration in the dissipative setting by focusing on models inspired by recent studies in celestial mechanics and cosmology.
}

\section{Introduction}

With the range of applications of contact geometry growing rapidly, geometric numerical integrators that preserve the contact structure have gained increasing attention~\cite{Bravetti2019,bravetti2020celestial,Pozs_r_2020,dediego2020discretecontact,vermeeren2019contact,zadra2020contact}.
Deferring to the above literature for detailed presentations of contact systems, their properties and many of their uses,
in this work we will present new applications of the contact geometric integrators introduced by the authors in~\cite{bravetti2020celestial,vermeeren2019contact,zadra2020contact} to two particular classes of examples inspired by celestial mechanics and cosmology.
\medskip

A contact manifold is a pair $(M,\xi)$ where $M$ is a $(2n+1)$-dimensional manifold and $\xi\subset{TM}$ is a contact structure, that is, a maximally non-integrable distribution of hyperplanes.
Locally, such distribution is given by the kernel of a one form $\eta$ satisfying $\eta \wedge (\dd \eta)^n \neq 0$ (see e.g.~\cite{geiges2008introduction} for more details). 
The 1-form $\eta$ is called the contact form. 
Darboux's theorem for contact manifolds states that for any point on $M$ there exist local coordinates $(q_1,\ldots,q_n,p_1,\ldots,p_n,s)$ such that the contact 1-form can be written as
$ \eta =  \dd s-\sum_i p_i \,\dd q_i . $
Moreover, given $\eta$, we can associate to any smooth function $ \mathcal{H}:M\rightarrow\mathbb{R}$ a contact Hamiltonian vector field $X_{\mathcal{H}}$, defined by
\begin{equation}\label{eq:XCH}
    \mathcal{L}_{X_{\mathcal{H}}} \eta = -R_{\eta}({\mathcal{H}}) \eta
\qquad \text{and} \qquad 
    \eta(X_{\mathcal{H}}) = - {\mathcal{H}},
\end{equation}
where $\mathcal{L}$ is the Lie derivative and $R_{\eta}$ is the Reeb vector field corresponding to $\eta$~\cite{geiges2008introduction}.
In canonical coordinates the flow of $X_{\mathcal{H}}$ is given by
\begin{equation}\label{eq:contactflow}
    \dot{q} = \der{\cH}{p}, \qquad
    \dot{p} = -\der{\cH}{q} - p \der{\cH}{s}, \qquad
    \dot{s} = p \der{\cH}{p} -{\mathcal{H}}.
\end{equation}
The flow of a contact Hamiltonian system preserves the contact structure, but it 
does not preserve the Hamiltonian:
\begin{equation}
\label{eq:hamevolution}
    \frac{\dd}{\dd t}{\mathcal{H}} = -{\mathcal{H}} \der{\cH}{s}.
\end{equation}
Using the this differential equation we can split the contact manifold in invariant parts for the Hamiltonian dynamics separated by the sub-manifold $\cH=0$, unique situation in which the Hamiltonian is conserved.

Contact Hamiltonian systems, like symplectic Hamiltonian systems, benefit from an associated variational principle, which is due to Herglotz: let $Q$ be an $n$--dimensional manifold with local coordinates $q^i$ and let $\cL: \R \times TQ \times \R \rightarrow \R$. For any given smooth curve $q:[0,T] \rightarrow Q$ we consider the initial value problem
\begin{equation}%
    \label{eq:action}
    \dot{s} = \cL(t,q(t),\dot{q}(t),s), \qquad s(0) = s_\mathrm{init}.
\end{equation}
Then the value $s(T)$ is a functional of the curve $q$. We say that $q$ is \emph{a critical curve} if $s(T)$ is invariant under infinitesimal variations of $q$ that vanish at the boundary of $[0,T]$.
It can be shown that critical curves for the Herglotz' variational principle are characterised by the following generalised Euler--Lagrange equations:
\begin{equation}%
    \label{eq:gEL}
    \der{\cL}{q^a} - \frac{\dd}{\dd t} \der{\cL}{\dot{q}^a} + \der{\cL}{s}\der{\cL}{\dot{q}^a} = 0.
\end{equation}
Furthermore, the corresponding flow consists of contact transformations with respect to the 1--form $\eta = \dd s - p_a \dd q^a$.
\medskip

This work is structured as follows: after a brief recap of the construction of contact integrators in Section 2, we showcase in Section 3 some interesting properties of the numerical integrators on some explicit examples. 
Finally, in Section 4 we present some considerations and ideas for future explorations.

\section{Contact integrators: theory}%
	\label{sec:contact}

    In this section we summarize the main results of~\cite{bravetti2020celestial,vermeeren2019contact} on the development of variational and Hamiltonian integrators for contact systems.

\subsection{Contact variational integrators (CVI)}%
	\label{sec:Lagrangian}
    There is a natural discretisation of Herglotz' variational principle~\cite{dediego2020discretecontact,vermeeren2019contact}.
    
    \begin{definition}[Discrete Herglotz' variational principle]
        Let $Q$ be an $n$--dimensional manifold with local coordinates $q^i$ and let $L: \R \times Q^2 \times \R^2 \rightarrow \R$. For any given discrete curve $q: \{0,\ldots,N\} \rightarrow Q$ we consider the initial value problem
        \(
        s_{k+1} = s_k + \tau L(k\tau,q_k,q_{k+1},s_k, s_{k+1})
        \), $s_0 = s_\mathrm{init}$.
        Then the value $s_N$ is a functional of the discrete curve $q$. We say that $q$ is 
        \emph{a critical curve} if
        \[ \der{s_N}{q_k} = 0 \qquad \forall k \in \{1,\ldots,N-1\}\,. \]
    \end{definition}

    From this, one can derive the discrete generalised Euler--Lagrange equations.
    As in the conventional discrete calculus of variations, they can be formulated as the equality of two formulas for the momentum~\cite{marsden2001discrete}.
    
    \begin{proposition}%
        \label{prop-dgEL}
    Let
    \begin{align*}
        p_k^- = \frac{ \displaystyle \partial_{q_k} L((k-1)\tau,q_{k-1},q_k,s_{k-1}, s_k) }{ \displaystyle 1 - \tau \partial_{s_k} L((k-1)\tau,q_{k-1},q_k,s_{k-1}, s_k) }\,, \quad
        p_k^+ = -\frac{ \displaystyle \partial_{q_k} L(k\tau,q_k,q_{k+1},s_k, s_{k+1}) }{ \displaystyle 1 + \tau \partial_{s_k} L(k\tau,q_k,q_{k+1},s_k, s_{k+1}) }\,.
        \end{align*}
    Then solutions to the discrete Herglotz variational principle are characterised by
    $p_k^- = p_k^+$.
    Furthermore, the map $(q_k,p_k,s_k) \mapsto (q_{k+1},p_{k+1},s_{k+1})$ induced by a critical discrete curve preserves the contact structure $\ker(\dd s - p_a \dd q^a)$.
    \end{proposition}
    
    Without loss of generality it is possible to take the discrete Lagrange function depending on only one instance of $s$: $L(k\tau,q_k,q_{k+1},s_k)$. Then the discrete generalised Euler--Lagrange equations read
    \begin{align*}
        \partial_{q_k} L(k\tau,q_k,q_{k+1},s_k)
        + \partial_{q_k} L((k-1)\tau,q_{k-1},q_k,s_{k-1}) \left( 1 + \tau \partial_{s_k} L(k\tau,q_k,q_{k+1},s_k) \right) = 0\,.
    \end{align*}

For a discrete Lagrangian of the form
\begin{equation}
    L(k\tau, q_k, q_{k+1}, s_k) = \frac12\left(\frac{q_{k+1} - q_{k}}\tau\right)^2 - \frac{V(q_{k}, k\tau) + V(q_{k+1}, (k+1)\tau)}2 - F(s_{k}, k\tau),
\end{equation}
which includes all the examples treated below, the CVI is explicit and takes the remarkably simple form
\begin{align}
    q_{k+1} &= q_{k} - \frac{\tau^2}{2} \partial_{q_k}V(q_{k}, k\tau) + p_{k}  \left(\tau - \frac{\tau^2}{2} \partial_{s_k} F(s_k, k\tau)\right) ,\\
s_{k+1} &= s_{k} + \tau L(k\tau, q_k, q_{k+1}, s_k) ,\\
p_{k+1} &= \frac{
    \left(1 - \frac\tau 2  \partial_{s_k}F(s_{k}, k\tau)\right)p_{k}
    - \frac{\tau}2 \left(\partial_{q_k}V(q_{k}, k\tau) + \partial_{q_{k+1}}V(q_{k+1}, (k+1)\tau)\right)
  }{1 + \frac{\tau}2 \partial_{s_{k+1}}F(s_{k+1}, (k+1)\tau)}.
\end{align}
{Higher order CVIs can be constructed with a Galerkin discretisation~\cite{bravetti2020celestial,dediego2020discretecontact}.}

\subsection{Contact Hamiltonian integrators (CHI)}%
    \label{sec:Hamiltonian}

In~\cite{bravetti2020celestial}, the authors
derive a contact analogue of the symplectic integrators introduced by Yoshida~\cite{yoshida1990construction} 
for separable Hamiltonians.

\begin{proposition}[\cite{bravetti2020celestial}]%
	\label{prop:2ndcontact}
    Let the contact Hamiltonian $\mathcal H$ be separable into the sum of $n$ functions $\phi_j(q,p,s)$, $j=1,\dots,n$.
    Assume that each of the
    vector fields $X_{\phi_j}$
    is exactly integrable and denote its flow by $\exp(t X_{\phi_j})$.
    Then the integrator
    \begin{equation}
        S_{2}(\tau)=e^{\frac{\tau}{2} X_{\phi_1}} e^{\frac{\tau}{2} X_{\phi_2}} \cdots e^{{\tau} X_{\phi_n}} \cdots e^{\frac{\tau}{2} X_{\phi_2}} e^{\frac{\tau}{2} X_{\phi_1}}, \label{secondinteghamil}
    \end{equation}
    is a second--order contact Hamiltonian integrator (CHI) for the flow of $\mathcal H$.
\end{proposition}

Note that exact integrability is not a necessary condition in general. The same construction gives in fact rise to \emph{contact compositional integrators} in a straightforward manner, but we will not digress further on this.

Starting from Proposition~\ref{prop:2ndcontact}, it is possible to construct CHIs of any even order: the construction is presented in detail in \cite{bravetti2020celestial},
where it was also shown how to use the \emph{modified Hamiltonian} in order to obtain error estimates for the integrator.

For a contact Hamiltonian of the form
\begin{equation}\label{eq:HamModl}
    \cH = \underbrace{\frac{p^2}{2}}_{C} + \underbrace{V(q,t)}_{B} + \underbrace{f(s,t)}_{A},
\end{equation}
one obtains for a time step $\tau$ the following discrete maps
\begin{align}
    A:  \begin{cases}
    q_i = q_{i-1} \\
    p_i = p_{i-1} \frac{f(s_{i},t_{i-1})}{f(s_{i-1},t_{i-1})}\\
    \int_{s_{i-1}}^{s_i} \frac{d \xi}{f(\xi,t_{i-1})} =- \tau
    \end{cases}
    B:  \begin{cases}
    q_i = q_{i-1} \\
    p_i = - V'(q_{i-1},t_{i-1}) \tau + p_{i-1}\\
    s_i = s_{i-1} - \frac{V(q_{i-1},t_{i-1})}{2} \tau
    \end{cases}
    C:  \begin{cases}
    q_i = p_{i-1} \tau + q_{i-1} \\
    p_i = p_{i-1}\\
    s_i = s_{i-1} + \frac{p^2_{i-1}}{2} \tau
    \end{cases}
\end{align}
The time is advanced by an extra map, $D: t_i = t_{i-1} + \tau$. 
Altogether, this can be seen in fluid dynamical terms as a streamline approximation of each piece of the Hamiltonian. All the examples in this paper are of the form \eqref{eq:HamModl} with an explicit map $A$.

\section{Contact integrators: applications}

Until recently, the main focus in the literature on contact mechanical systems has been on models with linear dependence on the action, see e.g.~\cite{bravetti2017contact,vermeeren2019contact}. 
In general, however, for a contact system to have non--trivial periodic trajectories in the contact space, $\frac{\partial \mathcal H}{\partial s}$ has to change sign.
This can be achieved either with a time--varying damping coefficient, as is the case in~\cite{Liu2018} and in our first example below, or by including a non--linear dependence on the action, as is done in~\cite{Huang2019} and in our second example.

Even though the contact oscillator that we use in the simulation is of purely theoretical interest and not associated to physical systems (that we know of), the reduction presented in \cite{sloan2018dynamical,Sloan2020A} shows that contact Hamiltonians with quadratic dependence on the action appear naturally in the study of the intrinsic dynamics of Friedman-Lemaitre-Robertson-Walker and Bianchi universes in cosmology.

The source code for all the simulations is provided in \cite{zenodoCodeGSI}.

\subsection{Perturbed Kepler problem}

Many relevant systems in celestial mechanics fall in the realm of Newtonian mechanics of systems with time--varying non--conservative forces \cite{bravetti2020celestial}.
Their equations of motion are the solution of the Newton equations
\begin{equation}%
	\label{eq:geneq}
    \ddot q + \frac{\partial V(q,t)}{\partial q}+f(t)\dot q=0\,.
\end{equation}
A direct computation, shows that \eqref{eq:geneq} coincide with the equations of motion of the contact Hamiltonian
\begin{equation}%
\label{eq:mainCH}
    \cH(p,q,s,t)=\sum_{a=1}^n\frac{p_a^2}{2}+V(q,t)+f(t)\,s.
\end{equation}
Since \eqref{eq:mainCH} is separable in the sense of Proposition~\ref{prop:2ndcontact}, one can directly apply the CVI and CHI to such Hamiltonian systems. 
Even in presence of singularities, like in the perturbed Kepler problem re-discussed here, contact integrators show a remarkable stability also for large time steps.

Here we consider a Kepler potential $V(q^a,t)=-\alpha/|q|$, $\alpha\in\R$, with an external periodic forcing $f(t)\dot q = \alpha \sin(\Omega t) \dot q$.
This choice of the perturbation is selected in order to have some sort of energy conservation on average, to emphasise the stability of the methods and their applicability to time-dependent problems.

It is well--known that, in presence of Keplerian potentials, Euclidean integrators become unstable for long integration times or large time steps: a Runge--Kutta integrator drifts toward the singularity and explodes in a rather short amount of iterations.

\begin{figure}[hbt]
    \centering
    \includegraphics[width=\onelinewidth]{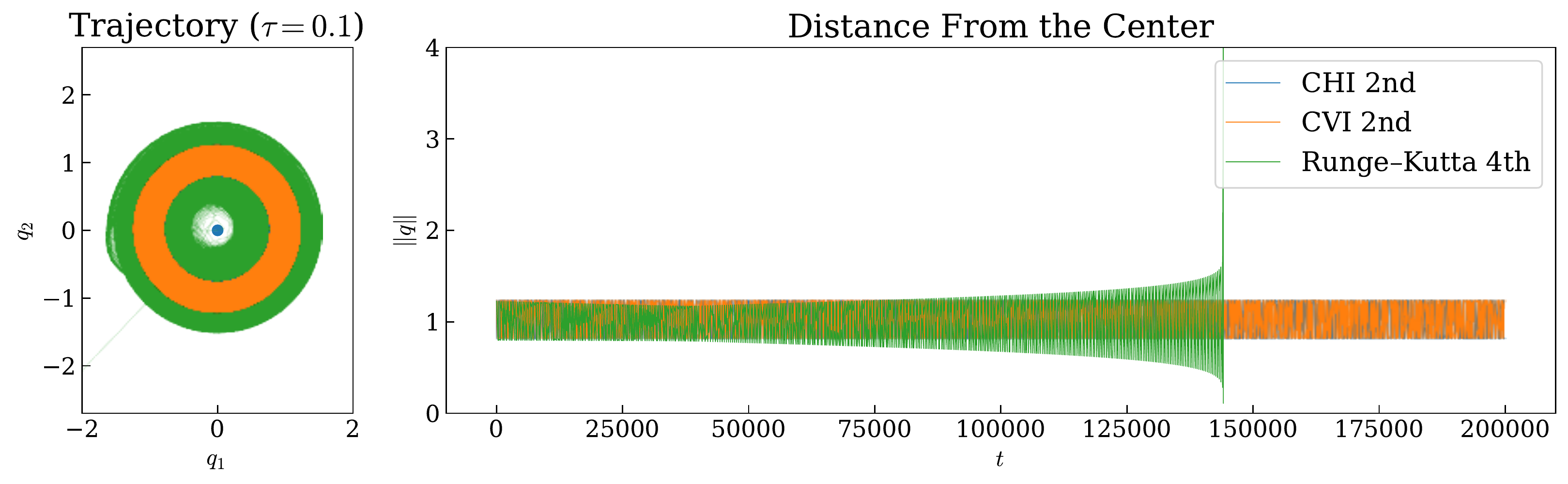}
    \caption{Long--time integration of a perturbed Kepler problem with time step $\tau\in\{0.1\}$.
    }%
    \label{fig:tpk}
\end{figure}

\begin{figure}[hbt]
    \centering
    \includegraphics[width=\onelinewidth]{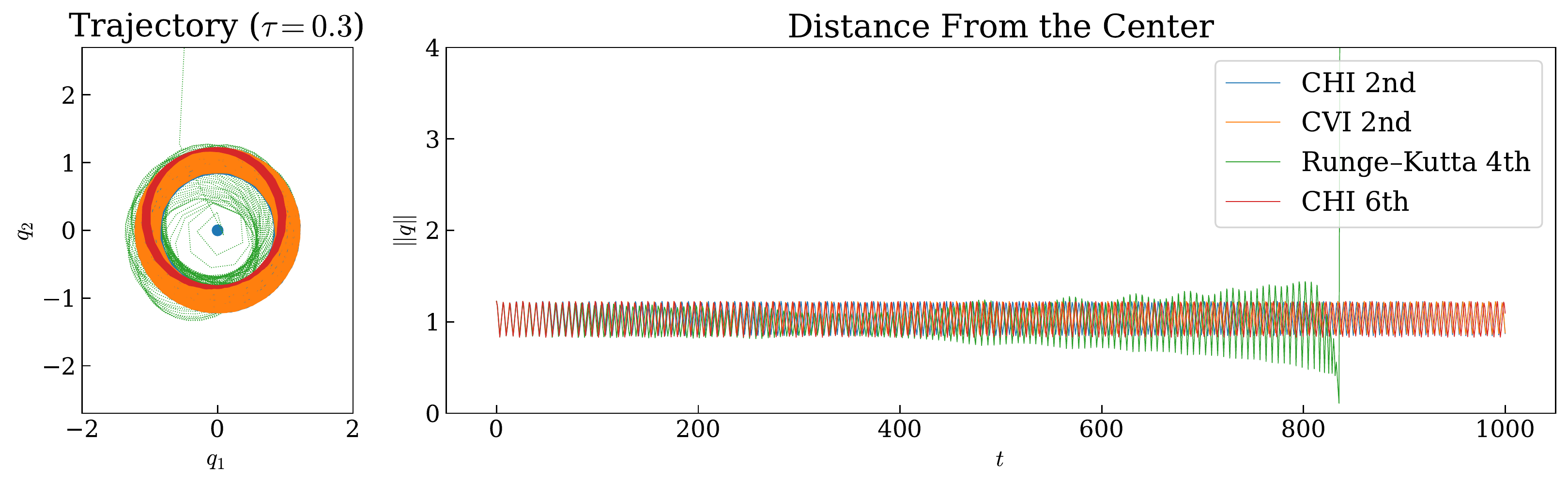}
    \caption{Integration of a perturbed Kepler problem with time step $\tau=0.3$. The plot is truncated after the first $1.000$ units of time to emphasize the blow-up of the Runge--Kutta 4 method.
    The contact integrators present no visible difference for the whole integration time.}%
	\label{fig:tpk0eprec}
\end{figure}

In Figure~\ref{fig:tpk}, the same perturbed Kepler problem with $\Omega = \pi$ and $\alpha=0.01$ is integrated over a very long time interval $[0,200.000]$ with the second--order contact integrators and a fixed--step fourth--order Runge--Kutta (RK4).
One can clearly observe the long-time stability of the contact method.
The price to pay for this 
is the introduction of an artificial precession of the trajectory. This is more evident in Figure~\ref{fig:tpk0eprec}: here the problem with $\Omega = \pi$ and $\alpha = 0.05$ is integrated with time step $\tau=0.3$ with the aforementioned integrators and with a 6th order CHI from~\cite{bravetti2020celestial}.

\subsection{Contact oscillator with quadratic action}\label{sec:qfric}

Motivated by the analysis in~\cite{Huang2019} and~\cite{sloan2018dynamical,Sloan2020A}, in this section we study
contact Hamiltonians of the form
$\cH(p,q,s)=\sum_{a=1}^n\frac{p_a^2}{2}+\,\gamma \frac{s^2}2 +V(q^a)$.
In particular, we consider the 1-dimensional quadratic contact harmonic oscillator, $V(q) = + \frac{q^2}{2} - C$, $\gamma , C>0$.
As shown by equation \eqref{eq:hamevolution}, the value of the contact Hamiltonian is not preserved unless its initial value is equal to zero \cite{Bravetti2017a}.
This generally defines an (hyper)surface in the contact manifold that separates two invariant basins for the evolution.
In the case at hand, the surface $\cH=0$ is an ellipsoid, or a sphere with radius $\sqrt{2 C}$ if $\gamma = 1$.
Furthermore, the quadratic contact oscillator presents two equilibrium points of different nature on $\cH=0$: the stable north pole $N = \left(0 , 0 , \sqrt{2 C \gamma^{-1}}\right)$ and the unstable south pole $S = \left(0 , 0 , -\sqrt{2 C\gamma^{-1}}\right)$.

In the case of geometric integrators, the explicit nature of the modified Hamiltonian allows to analyze this process and confirm that $\cH=0$ remains, for $\tau$ small enough, bounded and close to the original unperturbed surface.
For our particular example, a direct analysis of the fixed points of the integrator allows to control analytically the equilibrium points of the numerical map: these are simply shifted to $\left( 0, 0, \pm \frac12 \sqrt{8 C \gamma^{-1} + \tau^2  C^2}\right)$, and maintain their stability.

\begin{figure}[hbt]
    \centering
    \includegraphics[width=0.33\linewidth]{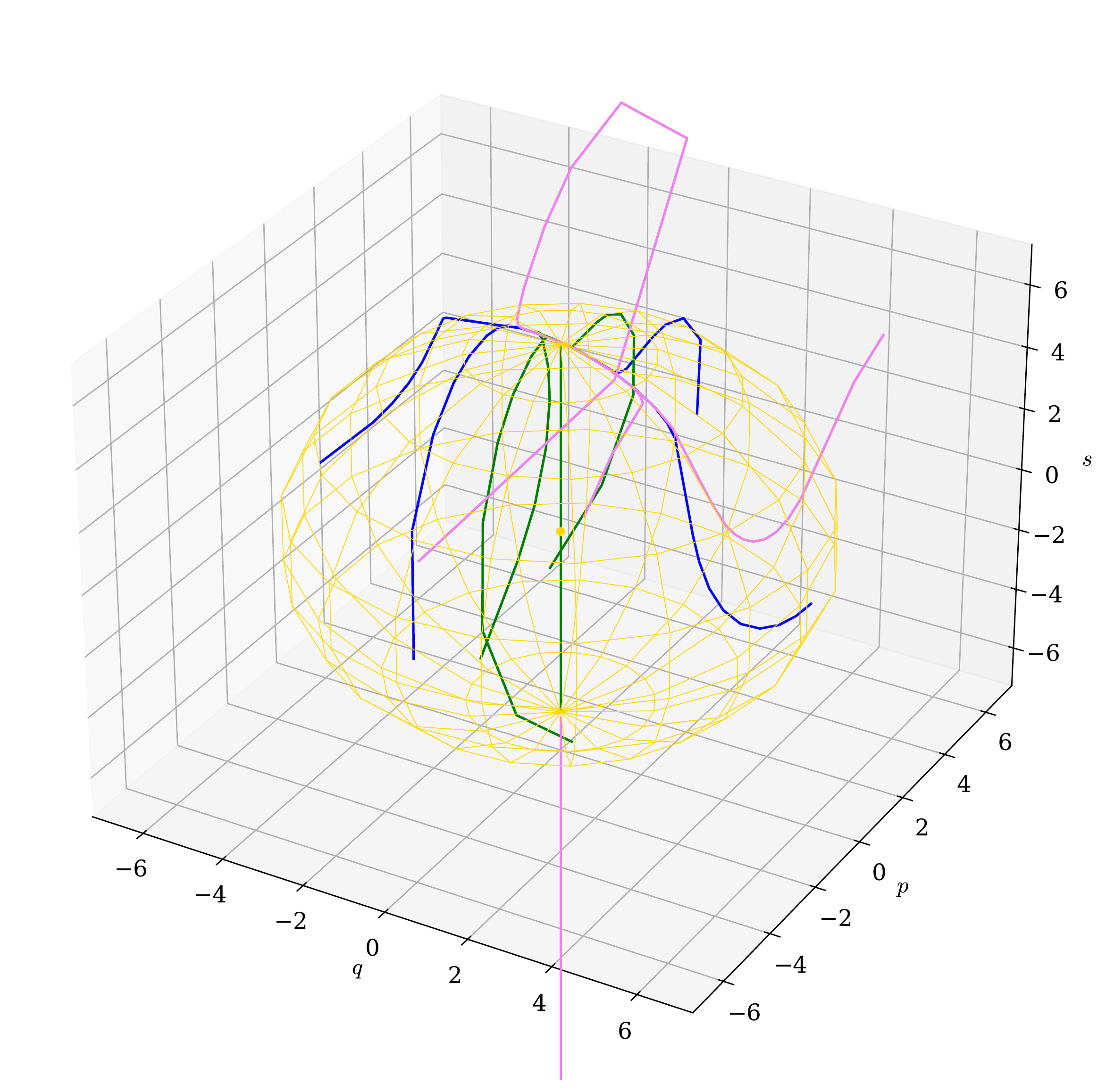}
    \includegraphics[width=0.33\linewidth]{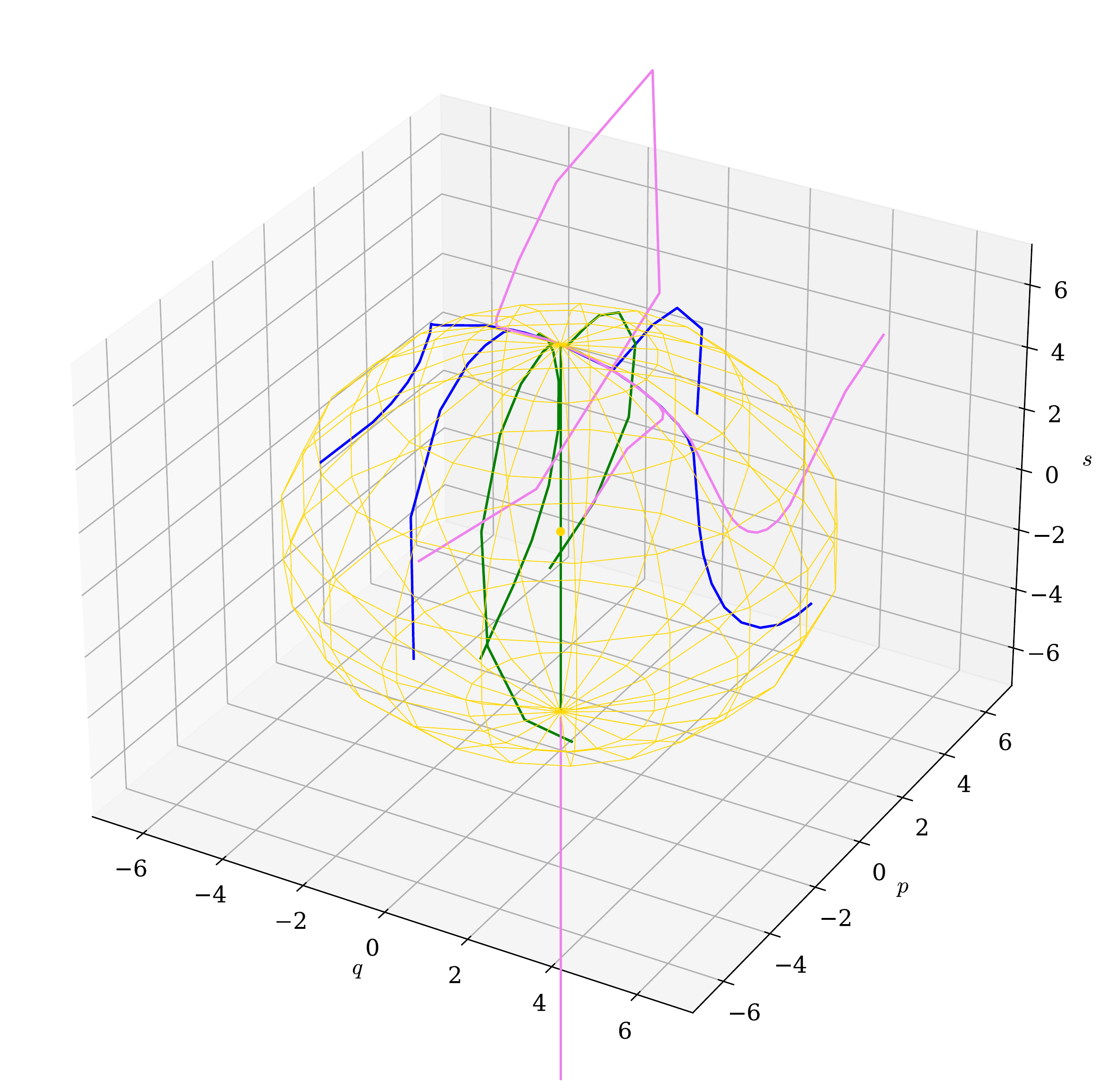}
    \caption{The yellow sphere of radius $6$ is the invariant surface $\cH=0$, for $\gamma=1$ and $C=18$.
    The trajectories are coloured according to their initial conditions: blue starts on the surface, purple outside and green inside.
    The half--line below $S$ is an unstable invariant submanifold of the system.
    Left is the CHI, right the CVI, both of 2nd order.}%
	\label{fig:sphere}
\end{figure}

In Figure~\ref{fig:sphere} we show the dynamics of trajectories starting on, outside and inside $\cH=0$.
The deformed invariant surface is so close to the sphere that they are practically indistinguishable even from a close analysis of the three dimensional dynamical plot.

\begin{figure}[hbt]
    \centering
    \includegraphics[width=\linewidth]{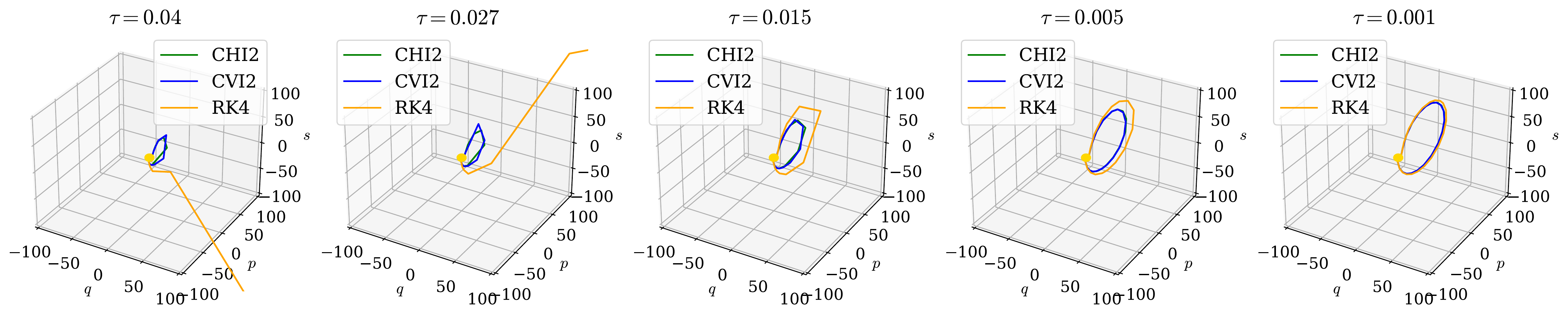}
    \caption{Trajectory with initial conditions $(q_0, p_0, s_0) = (0, -1, -7)$ integrated with different time steps. 
    Close to the unstable fixed point $(0,0,-6)$, the system becomes stiff: here the Runge--Kutta 4 integrator shows a higher degree of instability compared to the contact methods.}%
	\label{fig:sphereerr}
\end{figure}

The unstable south pole provides an excellent opportunity to compare the performance of our low--order integrators against the common Runge--Kutta 4 method:
initial conditions close to the unstable point are subject to fast accelerations away from the sphere, making the problem stiff.
One can clearly see this in Figure~\ref{fig:sphereerr}: for very small time steps we see that the trajectories are all converging towards the north pole, but as the time step gets larger we see the instability overcoming the Runge--Kutta integrator first and for larger $\tau$ also the CVI.
For this problem, in fact, the CHI displays a remarkable stability for larger times steps. Moreover, in terms of vector field evaluations, the low-order contact integrators are comparable to an explicit midpoint integrator \cite{zadra2020contact} and this reflects also in the  comparable running time (see Table~\ref{tbl:run}).

\begin{table}[hbt]
\centering
\begin{tabular}{@{}lll@{}}
\toprule
Integrator type (order)   & Mean time (from 10 runs) & Standard deviation \\ \midrule
CHI (2nd)                 & 0.0986                   & $\pm$ 0.0083         \\
CVI (2nd)                 & 0.0724                   & $\pm$ 0.0052         \\
Runge--Kutta (4th)        & 0.1375                   & $\pm$ 0.0075         \\
Midpoint (2nd)            & 0.0363                   & $\pm$ 0.0031         \\ \bottomrule
\end{tabular}
\label{tbl:run}
\caption{Integration of a contact oscillator with $\tau=0.1$ and $t\in[0,500]$.}
\end{table}

\section{Conclusions}%
	\label{sec:conclusions}
	
	In this manuscript we discussed some new directions for contact integrators. Even though from a physical perspective we presented basic examples, they show some of the generic advantages provided by contact integrators.
	In particular, they show the remarkable stability of our low--order integrators in comparison to standard higher--order methods.
	
	The study of contact integrators started in \cite{bravetti2020celestial,vermeeren2019contact,zadra2020contact} is still at its early stages. The recent work \cite{dediego2020discretecontact} has shown that they are underpinned by a beautiful geometric construction which appears to be closely related to non--holonomic and sub--Riemannian systems and which will require further investigation on its own.
	
    Furthermore, contact integrators for systems with non--linear dependence in the action might boost the analysis of new systems that can be of interest both for their dynamical structure and for their modelling capabilities, as it is already happening with their use in molecular dynamics~\cite{bravetti2016thermostat}, Monte Carlo algorithms~\cite{betancourt2014adiabatic} and relativistic  cosmology~\cite{sloan2018dynamical,Sloan2020A,gryb2021scale}.

\bibliographystyle{splncs04}
\bibliography{contact_int.bib}

\end{document}